\documentclass[reqno]{amsart}
\usepackage{amssymb}
\usepackage{amsmath}
\usepackage{amsfonts}
\usepackage{graphicx}

\newcommand{\vertiii}[1]{{\left\vert\kern-0.25ex\left\vert\kern-0.25ex\left\vert #1 
    \right\vert\kern-0.25ex\right\vert\kern-0.25ex\right\vert}}

\newtheorem{thm}{Theorem}[section]

\newtheorem{cor}[thm]{Corollary}
\newtheorem{lem}[thm]{Lemma}
\newtheorem{rem}[thm]{Remark}
\newtheorem{con}[thm]{Conjecture}

\newtheorem{note}[thm]{Note}
\newtheorem{que}[thm]{Question}

\newtheorem{ex}[thm]{Example}

\newtheorem{define}[thm]{Definition}
\numberwithin{equation}{section}
\newcommand{\RR}{\mathbb{R}}

\newcommand{\CC}{\mathbb{C}}

\newcommand{\tr}{\operatorname{tr~}}

 \newcommand{\MM}{\mathbb{M}}
 \newcommand{\HH}{\mathbb{H}}
 \newcommand{\PP}{\mathbb{P}}
\newcommand{\ie}{\text{,~i.e.,~}} 
 \newcommand{\BOX}{\hfill $\Box$}

\def\alert#1{\smallskip{\hskip\parindent\vrule\vbox{\advance\hsize-2\parindent\hrule\smallskip\parindent.4\parindent\narrower\noindent#1\smallskip\hrule}\vrule\hfill}\smallskip}

\newcommand{\mythe}[2]{
\begin{thm}
\ 
\label{#1}#2
\end{thm}}

\newcommand{\mycon}[2]{
\begin{con}
\ 
\label{#1} #2
\end{con}}

\newcommand{\mylem}[2]{
\begin{lem}
\ 
\label{#1} #2
\end{lem}}

\newcommand{\mycor}[2]{ 
\begin{cor}
\ 
\label{#1}#2
\end{cor}}

\newcommand{\myque}[2]{
\begin{que}
\ 
\label{#1}\rm #2
\end{que}}

\begin{document}
{\small{~~\\ \vspace{-2.5cm}\begin{flushright} \jobname\\ \today \end{flushright}\vspace{1.cm}}}

\title[On norm inequalities related to the geometric mean]{On norm inequalities related to the geometric mean}

\author[Freewan]{Shaima'a Freewan$^{1}$}
\address{$^{1}$Department of Mathematics, Yarmouk University, Irbid, Jordan}
\email{shyf725@gmail.com}
\author[Hayajneh]{Mostafa Hayajneh$^{2}$\ }
\address{$^{2}$Department of Mathematics, Yarmouk University, Irbid, Jordan}
\email{hayaj86@yahoo.com}


\begin{abstract}
Let $A_i$ and $B_i$ be positive definite matrices  for all $i=1,\cdots,m.$ It is shown that $$\left|\left|\sum_{i=1}^m(A_i^2\sharp B_i^2)^r\right|\right|_1\leq\left|\left|\left(\left(\sum_{i=1}^mA_i\right)^{\frac{pr}{_2}}\left(\sum_{i=1}^mB_i\right)^{pr}\left(\sum_{i=1}^mA_i\right)^{\frac{rp}{_2}}\right)^{\frac{1}{p}}\right|\right|_1,$$for all $p>0$ and for all $r\geq1.$

We conjecture this  inequality is also true for all unitarily invariant  norms. We give an affirmative answer to the case of  $m=2,$ $p\geq1$, $r\geq1$ and for all unitarily invariant norms. In other words,  it is shown that  $$\left|\left|\left|\left(A^{^2}\sharp B^{^2}\right)^{r}+\left(C^{^2}\sharp D^{^2}\right)^{r}\right|\right|\right|\leq \left|\left|\left|\left(\left(A+C\right)^{^\frac{rp}{_2}}\left(B+D\right)^{{rp}}\left(A+C\right)^{^\frac{rp}{_2}}\right)^{\frac{1}{_p}}\right|\right|\right|,$$for all  unitarly invariant norms, for all $p\geq1$ and for all $r\geq1$, where    $A,B,C,D$ are positive definite matrices. 
This gives an affirmative answer to the conjecture posed by  Dinh, Ahsani and Tam  in the case of  $m=2$.  The preceding inequalities directly lead to a recent result of Audenaert \cite{ANIFP}.
\end{abstract}

\maketitle

\section{Introduction}
Throughout this paper, we  refer to the set of all matrices of size $n$ and complex entries as $\MM_n(\CC)$.  The symbol $\HH_n$ indicates the set of all Hermitian matrices in $\MM_n(\CC)$.   If $A\in\HH_n$ is positive definite matrix\ie if $\langle x, Ax\rangle > 0$ for all $x\in\CC^n-\{0\}$, we use the notation $A>0$.  The symbol $\PP_n$ indicates the set of all  positive definite matrices in $\MM_n(\CC)$. 
In this paper, any unitarily invariant norm on the space $\MM_n(\CC)$ is denoted by the symbol  $\vertiii{.}.$ 
The Schatten $p$-norm is one of the most significant unitarily invariant norms. It is represented by the symbol $\left|\left|\cdot\right|\right|_p$, and is defined as 
\begin{eqnarray}
&&\left|\left|A\right|\right|_p=\left(\tr \left(|A|^{^p}\right)\right)^\frac{1}{_p} \text{ for all } 1\leq p\leq\infty.
\end{eqnarray}
 The trance norm, the Hibert-Schmidt norm and the spectral norm are represented by the $p$-norm for the values $p=1$, $p=2$ and $p=3$, respectively.
 Let  $A,B\in\PP_n.$ Then the matrix $A\sharp B$ is given by \begin{eqnarray}\label{A471}
A\sharp_\frac{1}{2} B=A^{^\frac{1}{_2}}\left( A^{^\frac{-1}{_2}}BA^{^\frac{-1}{_2}}\right)^{\frac{1}{_2}} A^{^\frac{1}{_2}}.
\end{eqnarray} This matrix is called the   geometric mean of $A$ and $B$. It was first mentioned in \cite{FCFS}, and is frequently abbreviated $A\sharp B$ in the literature. See  \cite{COCM}, \cite{PDM}, \cite{TRMOP} and \cite{FCFS}.  For all $t\in[0,1],$ the matrix $A\sharp_t B$ is given by \begin{eqnarray}\label{A4732}
A\sharp_t B=A^{^\frac{1}{_2}}\left( A^{^\frac{-1}{_2}}BA^{^\frac{-1}{_2}}\right)^{t} A^{^\frac{1}{_2}}.
\end{eqnarray} This matrix is called the $t$-geometric mean of $A$ and $B$.  

In \cite[Theorem 1.1]{AMS}, Bourin and Uchiyama  proved that \begin{eqnarray}\label{E1}
\left|\left|\left|f\left(A+B\right)\right|\right|\right|\leq\left|\left|\left|f\left(A\right)+f\left(B\right)\right|\right|\right|,
\end{eqnarray} for all positive semidefinite matrices $A,B\in\MM_n(\CC)$, for all nonnegative concave function $f:[0,\infty)\longrightarrow[0,\infty)$  and for all unitarily invariant norms. For nonnegative concave functions $f$ on $[0,\infty)$, the  inequality (\ref{E1}) is a noncommutative version of the well-known inequality \begin{eqnarray}\label{E2}
f(a + b)\leq f(a) + f(b), ~~ \forall a, b \geq0.
\end{eqnarray} In \cite{MSIAB}, Bourin asked the following question. 
\myque{E3}{If  $A,B\in\MM_n(\CC)$  are positive semidefinite matrices and $p,q>0,$ is it true that \begin{eqnarray*}
\left|\left|\left|A^{p+q}+B^{p+q}\right|\right|\right|\leq\left|\left|\left|(A^{p}+B^{p})(A^{q}+B^{q})\right|\right|\right|?
\end{eqnarray*}}
 In \cite{TIAAQ}, Hayajneh and Kittaneh provided an affirmative answer for the trace norm  and the Hilbert-Schmidt norm. In \cite{ANIFP}, Audenaert gave an affirmative answer to Question \ref{E3}  by proving  \begin{eqnarray}\label{E4}
\left|\left|\left|\sum_{i=1}^m\left(A_i B_i\right)\right|\right|\right|\leq\left|\left|\left|\sum_{i=1}^m\left(A_i^\frac{1}{_2}B_i^\frac{1}{_2}\right)^2\right|\right|\right|\leq\left|\left|\left|\left(\sum_{i=1}^mA_i\right)\left(\sum_{i=1}^mB_i\right)\right|\right|\right|,
\end{eqnarray} for all $A_i,B_i\in\PP_n,~i=1,\cdots,m$ such that $A_iB_i=B_iA_i$ and  all  unitarily invariant norms.
Also, this result confirms a conjecture of  Hayajneh and Kittaneh in \cite{TIAAQ}. In \cite{ROTRR}, Lin recently provided yet another proof of inequality (\ref{E4}). 

In \cite{TIAAQ2}, Hayajneh and Kittaneh generalized inequality (\ref{E4}) by proving  \begin{eqnarray}\label{E5}
\left|\left|\left|\sum_{i=1}^m\left(A_i B_i\right)\right|\right|\right|\leq\left|\left|\left|\sum_{i=1}^m\left(A_i^\frac{1}{_2}B_i^\frac{1}{_2}\right)^2\right|\right|\right|\leq\left|\left|\left|\left(\sum_{i=1}^mA_i\right)^\frac{1}{2}\left(\sum_{i=1}^mB_i\right)\left(\sum_{i=1}^mA_i\right)^\frac{1}{2}\right|\right|\right|,
\end{eqnarray} for all $A_i,B_i\in\PP_n,~i=1,\cdots,m$ such that $A_iB_i=B_iA_i$ and  all  unitarily invariant norms.
 
 In  \cite[Theorem~$3.1$]{GAI}, Dinh, Ahsani and Tam  proved the following remarkable inequality which is a non-commutative version of inequality (\ref{E5}),
 \begin{eqnarray}\label{A1}
\left|\left|\left|\sum_{i=1}^m\left(A_i\sharp B_i\right)^2\right|\right|\right|\leq\left|\left|\left|\left(\sum_{i=1}^mA_i\right)^\frac{1}{2}\left(\sum_{i=1}^mB_i\right)\left(\sum_{i=1}^mA_i\right)^\frac{1}{2}\right|\right|\right|,~~~
\end{eqnarray}  
for all $A_i,B_i\in\PP_n,$ $i=1,\cdots,m,$ and for all unitarily invariant norms. Dinh~recently gave a generalization in \cite{AIF}  by proving \begin{eqnarray}\label{A2}
\left|\left|\left|\sum_{i=1}^m\left(A_i\sharp_t B_i\right)^r\right|\right|\right|\leq\left|\left|\left|\left(\left(\sum_{i=1}^mA_i\right)^\frac{tpr}{2}\left(\sum_{i=1}^mB_i\right)^{(1-t)pr}\left(\sum_{i=1}^mA_i\right)^\frac{tpr}{2}\right)^\frac{1}{p}\right|\right|\right|,\hspace{0.3cm}
\end{eqnarray}
for all  $t\in[0,1],$  for all $r\geq1,$  for all $ p>0$,  for all $A_i,B_i\in\PP_n,$ $i=1,\cdots,m,$ and  for all unitarily invariant norms. In \cite[page 787]{GAI}, Dinh, Ahsani and Tam  proposed the following conjecture.
\mycon{CON5}{Let $A_i,B_i\in\PP_n$ for all  $i=1,\cdots,m.$ Then for  all unitarily invariant norms \begin{eqnarray}\label{A492}
\left|\left|\left|\sum_{i=1}^m\left(A_i^2\sharp B_i^2\right)\right|\right|\right|\leq\left|\left|\left|\left(\sum_{i=1}^mA_i\right)^\frac{1}{2}\left(\sum_{i=1}^mB_i\right)\left(\sum_{i=1}^mA_i\right)^\frac{1}{2}\right|\right|\right|.
\end{eqnarray}}
The authors of the same paper \cite[Corollary 3.3]{GAI}  proved Conjecture \ref{CON5} for the case of the trace norm $\|.\|_1.$ 
In this paper, we present the following more general conjecture.
\mycon{CON2}{
Let $A_i,B_i\in\PP_n$ for all $i=1,\cdots,m.$ Then for  all unitarily invariant norms $$\left|\left|\left|\sum_{i=1}^m(A_i^2\sharp_t B_i^2)^r\right|\right|\right|\leq\left|\left|\left|\left(\left(\sum_{i=1}^mA_i\right)^{\frac{prt}{_2}}\left(\sum_{i=1}^mB_i\right)^{2(1-t)pr}\left(\sum_{i=1}^mA_i\right)^{\frac{trp}{_2}}\right)^{\frac{1}{p}}\right|\right|\right|,$$for all  $t\in[0,1],$  for all  $r\geq1$ and for all $p>0.$}
 In Section \ref{S22},  we provide an affirmative answer to  Conjecture \ref{CON2} for the case of $m=2,$ $p\geq1$, $r\geq1$, $t=\frac{1}{2}$ and for all unitarily invariant norms.  In other words, if  $A,B,C,D\in\PP_n$, then for all $p\geq1$, for all $r\geq1$  and for  all unitarily invariant norms \begin{eqnarray}\label{A3}
\left|\left|\left|\left(A^{^2}\sharp B^{^2}\right)^{r}+\left(C^{^2}\sharp D^{^2}\right)^{r}\right|\right|\right|\leq \left|\left|\left|\left(\left(A+C\right)^{^\frac{rp}{_2}}\left(B+D\right)^{{rp}}\left(A+C\right)^{^\frac{rp}{_2}}\right)^{\frac{1}{_p}}\right|\right|\right|.
\end{eqnarray} In particular, this result settles  Conjecture \ref{CON5}  in the case of  $m=2$.
In Section \ref{S33}, we also give an affirmative answer to Conjecture \ref{CON2} for the case of the trace norm $\|.\|_1$. In other words, we prove $$\left|\left|\sum_{i=1}^m(A_i^2\sharp B_i^2)^r\right|\right|_1\leq\left|\left|\left(\left(\sum_{i=1}^mA_i\right)^{\frac{pr}{_2}}\left(\sum_{i=1}^mB_i\right)^{pr}\left(\sum_{i=1}^mA_i\right)^{\frac{rp}{_2}}\right)^{\frac{1}{p}}\right|\right|_1,$$
for all  $A_i,B_i\in\PP_n,$ $i=1,\cdots,m,$ for all $p>0$ and  for all $r\geq1.$ 
 In Section \ref{S44},  we apply our approach to give a new proof of inequality (\ref{A2})  for the case of $p=m=2,$ $t=\frac{1}{2}$ and  $r\geq1\ie$the following inequality,
 \begin{eqnarray}\label{A9}
\left|\left|\left|(A\sharp B)^r+(C\sharp D)^r\right|\right|\right|\leq\left|\left|\left|\left(\left(A+C\right)^\frac{r}{_2}\left(B+D\right)^r\left(A+C\right)^{^\frac{r}{_2}}\right)^{^\frac{1}{_2}}\right|\right|\right|.
\end{eqnarray}

\section{Preliminary}\label{S11}
We begin with some basic, well-known facts that will be utilized to prove our major results.

Let $A\in\MM_n(\CC)$.  We denote the {\it absolute value} of $A$ by $|A|=(A^*A)^{^\frac{1}{_2}}.$  The  singular values of $A$ are the $n$-nonnegative numbers $s_1(A),s_2(A), \cdots,s_n(A)$ such that   for all $i=1,2,\cdots,n$\begin{eqnarray*}
s_i(A)=\lambda_i(|A|)=\lambda_i^{^\frac{1}{_2}}(A^*A)=\lambda_i^{^\frac{1}{_2}}(AA^*),
\end{eqnarray*} where $\lambda_1(A), \cdots, \lambda_n(A)$ are the eigenvalues  of $A.$  Let  $s_1(A), \cdots, s_n(A)$ be arranged in such a way that $s_1(A)\geq \cdots\geq s_n(A).$ Then we will write $s(A)$ for the $n$-vector with components $s_1(A), \cdots, s_n(A)$\ie $$s(A)=(s_1(A),s_2(A),\cdots,s_n(A)).$$  Let $A\in\HH_n$ with eigenvalues $\lambda_1(A), \cdots, \lambda_n(A)$ arranged in such a way that $\lambda_1(A)\geq\cdots\geq\lambda_n(A).$ Then we will write $\lambda(A)$ for the $n$-vector with components $\lambda_1(A),  \cdots,\\\lambda_n(A)$\ie
$$\lambda(A)=(\lambda_1(A),\lambda_2(A),\cdots,\lambda_n(A)).$$

 Now, we start with the concept of majorisation.
Let $x = (x_1,\cdots, x_n)$ and $y = (y_1,\cdots, y_n)$ be two $n$-tuples of real numbers. Let
$x^{\downarrow}_1\geq x^{\downarrow}_2\geq\cdots\geq x^{\downarrow}_n$ be the decreasing rearrangement of $x_1, x_2,\cdots,x_n$. If for all $1 \leq k \leq n$ \begin{eqnarray}\label{A382}
\sum_{j=1}^kx^{\downarrow}_j\leq\sum_{j=1}^ky^{\downarrow}_j,
\end{eqnarray}then we say that $x$ is weakly  majorised by $y,$ and write this as $x\prec_{w}y$.  If, in addition to (\ref{A382}) we also have \begin{eqnarray}\label{A383}
\sum_{j=1}^nx^{\downarrow}_j=\sum_{j=1}^ny^{\downarrow}_j,
\end{eqnarray}then we say that $x$ is  majorised by $y,$ and write this as $x\prec y$.
Let $x = (x_1,\cdots, x_n),~y = (y_1,\cdots, y_n)$ be two $n$-tuples of nonnegative numbers. Let
$x^{\downarrow}_1\geq x^{\downarrow}_2\geq\cdots\geq x^{\downarrow}_n$ be the decreasing rearrangement of $x_1, x_2,\cdots,x_n$. If for all $1 \leq k \leq n$ \begin{eqnarray}\label{A380}
\prod_{j=1}^kx^{\downarrow}_j\leq\prod_{j=1}^ky^{\downarrow}_j,
\end{eqnarray} then we say that $x$ is weakly log majorised by $y.$ We write this as $x\prec_{w\log}y$.  If, in addition to (\ref{A380}) we also have \begin{eqnarray}\label{A381}
\prod_{j=1}^nx^{\downarrow}_j=\prod_{j=1}^ny^{\downarrow}_j,
\end{eqnarray}then we say that $x$ is log majorised by $y.$ We write this as $x\prec_{\log}y$.

Lemma \ref{A409} and Lemma \ref{A341}  give  some properties of majorisation inequalities.
\mylem{A409}{Let $x ,y$ be two $n$-tuples of nonnegative numbers. Then
$$x\prec_{\log} y\Longrightarrow x\prec_{w\log} y\Longrightarrow x\prec_{w}y.$$}
{\it Proof.~}See  \cite[pages 345]{MT}.\BOX

\mylem{A341}{Let $x = (x_1,\cdots, x_n),~y = (y_1,\cdots, y_n)$ be two $n$-tuples of nonnegative numbers. Then for all $r\geq1$ \begin{eqnarray*}
x\prec_wy\Longrightarrow x^r\prec_wy^r\ie(x_1^r,\cdots,x_n^r)\prec_w(y_1^r,\cdots,y_n^r).
\end{eqnarray*}}
{\it Proof.~}See \cite[page 42]{MA} and \cite[page 342]{MT}.\BOX

The following two lemmas give some properties of the $t$-geometric mean. 
\mylem{C2}{If $A,B\in\PP_n$,  then there exists a unitary~ $U\in~\MM_n(\CC)$ such that 
$A\sharp B = A^{^\frac{1}{_2}}UB^{^\frac{1}{_2}}.$}
{\it Proof.~}See \cite[page 108]{PDM}.\BOX

The next theorem is devoted to a well-known   log-majorisation inequality.
\mylem{A403}{(Ando-Hiai)~Let $A, B\in\PP_n$ and $t\in[0, 1]$. Then 
\begin{eqnarray*}
\lambda\left(A^r\sharp_tB^r\right)\prec_{\log}\lambda\left(\left(A\sharp_tB\right)^r\right),\hspace{0.5cm}\forall r\geq1.
\end{eqnarray*}}
{\it Proof.~}See \cite[pages 119–120]{LMA}.\BOX

Lemma \ref{A294} and Lemma \ref{A298} are devoted to present two powerful weakly majorized inequalities.
\mylem{A294}{(Weyl's Majorant Theorem)~Let $A\in\MM_n(\CC)$   with eigenvalues $\lambda_1(A),\cdots,\lambda_n(A)$ arranged in
such a way that $|\lambda_1(A)|\geq \cdots\geq|\lambda_n(A)|.$ Then for all $p\geq0$, we have \begin{eqnarray*}\left|\lambda(A)\right|^p\prec_w\lambda\left(\left|A\right|^p\right).
\end{eqnarray*}
In  other words,\begin{eqnarray}
\sum_{i=1}^k|\lambda_i(A)|^p\leq\sum_{i=1}^k\lambda_i^p(|A|),   ~~\forall  k=1,2,\cdots,n.
\end{eqnarray}}
{\it Proof.~}See \cite[page 42]{MA}.\BOX
\mylem{A298}{Let $A, B\in\MM_n(\CC)$. Then for all $p>0$, we have 
\begin{eqnarray*}
s^{p}(AB)\prec_w~s^{p}(A)s^{p}(B).
\end{eqnarray*}
 In other words,
\begin{eqnarray*}
\sum_{i=1}^k s_i^{p}(AB)\leq\sum_{i=1}^k s_i^{p}(A)s_i^{p}(B),  \hspace{0.5cm}\forall k=1,2,\cdots,n.
\end{eqnarray*}  }
{\it Proof.~}See \cite[page 94]{MA}.\BOX

Let us now present three useful lemmas regarding unitarily invariant norm inequalities.
\mylem{C32}{(Fan Dominance Theorem)~Let $A,B\in\MM_n(\CC)$. Then\begin{eqnarray*}
s(A)\prec_ws(B)\Longleftrightarrow|||A||| \leq|||B||| \hspace{0.2cm}\text{ for all unitarily invariant norms.}
\end{eqnarray*}}
{\it Proof.~}See \cite[page 375]{MT} and \cite[page 93]{MA}.\BOX

\mylem{A306}{Let $A, B\in\MM_n(\CC)$  such that the product
$AB$ is normal. Then~for every unitarily invariant norms, we have $|||AB||| \leq|||BA|||.$}
{\it Proof.~}See \cite[page 253]{MA}.\BOX

\mylem{A429}{Let $A,B\in\PP_n$ and $r\geq1.$  Then for all unitarily invariant norms $$\left|\left|\left|A\right|\right|\right|\leq\left|\left|\left|B\right|\right|\right|\Longrightarrow\left|\left|\left|A^{^r}\right|\right|\right|\leq \left|\left|\left|B^{^r}\right|\right|\right|.$$}
{\it Proof.~}Note that
\begin{eqnarray*}
\left|\left|\left|A\right|\right|\right|\leq\left|\left|\left|B\right|\right|\right|
&\Longrightarrow&s(A)\prec_{w}s(B)\hspace{0.5cm}\text{(by Lemma \ref{C32})}\\
&\Longrightarrow&s^r(A)\prec_{w}s^r(B)\hspace{0.5cm}\text{(by Lemma \ref{A341})}\\
&\Longrightarrow&s\left(A^{^r}\right)\prec_{w}s\left(B^{^r}\right)\\
&\Longrightarrow&\left|\left|\left|A^{^r}\right|\right|\right|\leq\left|\left|\left|B^{^r}\right|\right|\right|.\hspace{0.5cm}\text{(by Lemma \ref{C32})}
\end{eqnarray*}
This completes the proof.\BOX

Our next theorem connects the concepts of convex functions and  unitarily~invariant norms.
\mylem{A389}{Let $A_i\geq 0$  for all $i=1,\cdots,m.$ Let $f:[0,\infty)\longrightarrow[0,\infty)$ be a convex function with $f(0)=0$.
Then for all unitarily invariant norms, we have
 \begin{eqnarray*}
 \left|\left|\left|\sum_{i=1}^mf(A_i)\right|\right|\right|\leq\left|\left|\left|f\left(\sum_{i=1}^mA_i\right)\right|\right|\right|.
\end{eqnarray*}}
{\it Proof.~}See \cite[Theorem 1.2]{AMS}.\BOX

Let $E=\left[\begin{array}{ccc}
 A & 0\\
0 & B
\end{array}\right]\in\MM_{2n}(\CC),$ in which $A,B\in\MM_{n}(\CC)$.  Then we say that the matrix $E$ is the {\it direct sum } of $A$ and $B$, and write it as $E=A\oplus B$. 

\mylem{A299}{Let $A,B\in\MM_n(\CC)$. Then for all unitarily invariant norms on $\MM_{2n}(\CC)$, we have\begin{eqnarray*}
\left|\left|\left|\left[\begin{array}{lll}
A&0\\
0&B
\end{array}\right]\right|\right|\right|\leq \left|\left|\left|\left[\begin{array}{lll}
|A|+|B|&0\\
0&0
\end{array}\right]\right|\right|\right|.
\end{eqnarray*}
 In other words,  for all $k=1,2,\cdots,2n,$ we have
\begin{eqnarray*}
\sum_{i=1}^k s_i\left(\left[\begin{array}{lll}
A&0\\
0&B
\end{array}\right]\right)\leq\sum_{i=1}^k s_i\left(\left[\begin{array}{lll}
|A|+|B|&0\\
0&0
\end{array}\right]\right).
\end{eqnarray*}}
{\it Proof.~}See \cite[page 97]{MA} and Lemma  \ref{C32}.\BOX

The next theorem introduces two  weakly log majorised  inequalities.
\mylem{A416}{If  $A,B\in\PP_n,$  then \begin{enumerate}
\item$s\left(A^{^r}B^{^r}A^{^r}\right)\prec_{w\log}s\left(\left(ABA\right)^{^r}\right)~\text{for  all }r\in~[0,1].$
\item$s\left(\left(ABA\right)^{^r}\right)\prec_{w\log}s\left(A^{^r}B^{^r}A^{^r}\right)~\text{for  all }r\geq1.$
\end{enumerate}}
{\it Proof.~}The inequality (1) was proved in \cite[page 258]{MA}. The inequality (2) is proved in exactly the same way as proof of the inequality  (1).  \BOX
  
 The following theorem  is a generalization of  Lemma \ref{A416}.

\mylem{A417}{Let  $A,B\in\PP_n$  and let $p>0$. Then 
\begin{enumerate}
\item $s\left(\left(ABA\right)^{rp}\right)\prec_{w\log}s\left(\left(A^{^r}B^{^r}A^{^r}\right)^p\right)$\text{  for all }$r\geq1$
\item $s\left(\left(A^{^r}B^{^r}A^{^r}\right)^p\right)\prec_{w\log}s\left(\left(ABA\right)^{rp}\right)$\text{ for all }$r\in[0,1]$.
\end{enumerate}}
{\it Proof.~}We are to prove (1).~Let~$A,B\in\PP_n$  and let $p>0$ and $r\geq1$.    Using (2) in Lemma \ref{A416},  we have $$\displaystyle\prod_{i=1}^ks_i\left(\left(ABA\right)^{^r}\right)\leq\displaystyle{\prod_{i=1}^ks_i}\left(A^{^r}B^{^r}A^{^r}\right)~\text{for all }k=1,\cdots,n.$$  
Now, we have \begin{eqnarray*}
&&\displaystyle\prod_{i=1}^ks_i\left(\left(ABA\right)^{^r}\right)\leq\displaystyle{\prod_{i=1}^ks_i}\left(A^{^r}B^{^r}A^{^r}\right)~\text{for all }k=1,\cdots,n\\
&\Longrightarrow&\displaystyle\prod_{i=1}^ks_i^p\left(\left(ABA\right)^{^r}\right)\leq\displaystyle\prod_{i=1}^ks_i^p\left(A^{^r}B^{^r}A^{^r}\right)~\text{for all }k=1,\cdots,n\\
&\Longrightarrow&s^p\left(\left(ABA\right)^{^r}\right)\prec_{w\log}s^p\left(A^{^r}B^{^r}A^{^r}\right)\\
&\Longrightarrow&s\left(\left(ABA\right)^{^{pr}}\right)\prec_{w\log}s\left(\left(A^{^r}B^{^r}A^{^r}\right)^p\right).
\end{eqnarray*}
So  $s\left(\left(ABA\right)^{^{rp}}\right)\prec_{w\log}s\left(\left(A^{^r}B^{^r}A^{^r}\right)^p\right).$
Similarly, we can prove part (2).  This completes the proof.\BOX

\mylem{A295}{Let $A\in\MM_n(\CC)$ with singular values $s_1(A), s_2(A), \cdots, s_n(A)$. Then  \begin{eqnarray*}
s_i(UAW)=s_i(A),
\end{eqnarray*} for all $i=1,2,\cdots,n$ and for  all unitaries $U,W\in~\MM_n(\CC).$}
{\it Proof.~}See \cite[page 91]{MA}.\BOX 
\section{Main Results}\label{S22}
We start with the following three lemmas and then state and prove our main theorem.
\mylem{C53}{Let $A,B,C,D\in\PP_n.$ Then  for all $i=1,2,\cdots,2n,$ we have  \begin{eqnarray*}
s_i^2\left(\left[\begin{array}{lll}
B^{^\frac{1}{_2}}&0\\
D^{^\frac{1}{_2}}&0
\end{array}\right]\left[\begin{array}{lll}
A^{^\frac{1}{_2}}&C^{^\frac{1}{_2}}\\
0&0
\end{array}\right]\right)=s_i\left(\left[\begin{array}{lll}
\left(A+C\right)^{^\frac{1}{_2}}\left(B+D\right)\left(A+C\right)^{^\frac{1}{_2}}&0\\
0&0
\end{array}\right]\right).
\end{eqnarray*}}
{\it Proof.~}Let $A,B,C,D\in\PP_n$ and  let  $i=1,2,\cdots,2n.$ We  are to show that \begin{eqnarray*}
s_i^2\left(\left[\begin{array}{lll}
B^{^\frac{1}{_2}}&0\\
D^{^\frac{1}{_2}}&0
\end{array}\right]\left[\begin{array}{lll}
A^{^\frac{1}{_2}}&C^{^\frac{1}{_2}}\\
0&0
\end{array}\right]\right)=s_i\left(\left[\begin{array}{lll}
\left(A+C\right)^{^\frac{1}{_2}}\left(B+D\right)\left(A+C\right)^{^\frac{1}{_2}}&0\\
0&0
\end{array}\right]\right).
\end{eqnarray*} To see this, note that
\begin{eqnarray*}
&&s_i^2\left(\left[\begin{array}{lll}
B^{^\frac{1}{_2}}&0\\
D^{^\frac{1}{_2}}&0
\end{array}\right]\left[\begin{array}{lll}
A^{^\frac{1}{_2}}&C^{^\frac{1}{_2}}\\
0&0
\end{array}\right]\right)\\
&=&\lambda_i\left(\left(\left[\begin{array}{lll}
B^{^\frac{1}{_2}}&0\\
D^{^\frac{1}{_2}}&0
\end{array}\right]\left[\begin{array}{lll}
A^{^\frac{1}{_2}}&C^{^\frac{1}{_2}}\\
0&0
\end{array}\right]\right)^*\left(\left[\begin{array}{lll}
B^{^\frac{1}{_2}}&0\\
D^{^\frac{1}{_2}}&0
\end{array}\right]\left[\begin{array}{lll}
A^{^\frac{1}{_2}}&C^{^\frac{1}{_2}}\\
0&0
\end{array}\right]\right)\right)\\
&=&\lambda_i\left(\left[\begin{array}{lll}
A^{^\frac{1}{_2}}&0\\
C^{^\frac{1}{_2}}&0
\end{array}\right]\left[\begin{array}{lll}
B^{^\frac{1}{_2}}&D^{^\frac{1}{_2}}\\
0&0
\end{array}\right]\left[\begin{array}{lll}
B^{^\frac{1}{_2}}&0\\
D^{^\frac{1}{_2}}&0
\end{array}\right]\left[\begin{array}{lll}
A^{^\frac{1}{_2}}&C^{^\frac{1}{_2}}\\
0&0
\end{array}\right]\right)\\
&=&\lambda_i\left(\left[\begin{array}{lll}
A^{^\frac{1}{_2}}&C^{^\frac{1}{_2}}\\
0&0
\end{array}\right]\left[\begin{array}{lll}
A^{^\frac{1}{_2}}&0\\
C^{^\frac{1}{_2}}&0
\end{array}\right]\left[\begin{array}{lll}
B^{^\frac{1}{_2}}&D^{^\frac{1}{_2}}\\
0&0
\end{array}\right]\left[\begin{array}{lll}
B^{^\frac{1}{_2}}&0\\
D^{^\frac{1}{_2}}&0
\end{array}\right]\right)\\
&=&\lambda_i\left(\left[\begin{array}{lll}
A+C&0\\
0&0
\end{array}\right]\left[\begin{array}{lll}
B+D&0\\
0&0
\end{array}\right]\right)\\
&=&\lambda_i\left(\left[\begin{array}{lll}
\left(A+C\right)^{^\frac{1}{_2}}&0\\
0&0
\end{array}\right]\left[\begin{array}{lll}
B+D&0\\
0&0
\end{array}\right]\left[\begin{array}{lll}
\left(A+C\right)^{^\frac{1}{_2}}&0\\
0&0
\end{array}\right]\right)\\
&=&s_i\left(\left[\begin{array}{lll}
\left(A+C\right)^{^\frac{1}{_2}}\left(B+D\right)\left(A+C\right)^{^\frac{1}{_2}}&0\\
0&0
\end{array}\right]\right).\hspace{9cm}
\end{eqnarray*}

This shows that  \begin{eqnarray*}
s_i^2\left(\left[\begin{array}{lll}
B^{^\frac{1}{_2}}&0\\
D^{^\frac{1}{_2}}&0
\end{array}\right]\left[\begin{array}{lll}
A^{^\frac{1}{_2}}&C^{^\frac{1}{_2}}\\
0&0
\end{array}\right]\right)=s_i\left(\left[\begin{array}{lll}
\left(A+C\right)^{^\frac{1}{_2}}\left(B+D\right)\left(A+C\right)^{^\frac{1}{_2}}&0\\
0&0
\end{array}\right]\right),
\end{eqnarray*}for all $i=1,2,\cdots,2n.$  This completes the proof.\BOX

\mylem{C51}{Let $A,B,C,D\in\PP_n$ and  let $r>0$ and $p>0$ such that  $pr\geq1.$  Then  for all $k=1,2,\cdots,2n,$ we have  \begin{eqnarray*}
\sum_{i=1}^{k}s_i^{2r}\left(\left[\scalebox{0.95}{$\begin{array}{lll}
B^{^\frac{1}{_2}}&0\\
D^{^\frac{1}{_2}}&0
\end{array}$}\right]\left[\begin{array}{lll}
A^{^\frac{1}{_2}}&C^{^\frac{1}{_2}}\\
0&0
\end{array}\right]\right)\leq\sum_{i=1}^{k}s_i\left(\left[\scalebox{0.88}{$\begin{array}{lll}
\left(\left(A+C\right)^{^\frac{rp}{_2}}\left(B+D\right)^{{rp}}\left(A+C\right)^{^\frac{rp}{_2}}\right)^{\frac{1}{_p}}&0\\
~~0&0
\end{array}$}\right]\right).
\end{eqnarray*}}
{\it Proof.~}Let $A,B,C,D\in\PP_n$ and  let  $k=1,2,\cdots,2n.$   Let $r>0$ and $p>0$ such that  $pr\geq1.$  We are to prove that  \begin{eqnarray*}
\sum_{i=1}^{k}s_i^{2r}\left(\left[\scalebox{0.95}{$\begin{array}{lll}
B^{^\frac{1}{_2}}&0\\
D^{^\frac{1}{_2}}&0
\end{array}$}\right]\left[\begin{array}{lll}
A^{^\frac{1}{_2}}&C^{^\frac{1}{_2}}\\
0&0
\end{array}\right]\right)\leq\sum_{i=1}^{k}s_i\left(\left[\scalebox{0.88}{$\begin{array}{lll}
\left(\left(A+C\right)^{^\frac{rp}{_2}}\left(B+D\right)^{{rp}}\left(A+C\right)^{^\frac{rp}{_2}}\right)^{\frac{1}{_p}}&0\\
~~0&0
\end{array}$}\right]\right).
\end{eqnarray*} From Lemma \ref{C53}, we get
\begin{eqnarray*}
s_i^2\left(\left[\begin{array}{lll}
B^{^\frac{1}{_2}}&0\\
D^{^\frac{1}{_2}}&0
\end{array}\right]\left[\begin{array}{lll}
A^{^\frac{1}{_2}}&C^{^\frac{1}{_2}}\\
0&0
\end{array}\right]\right)=s_i\left(\left[\begin{array}{lll}
\left(A+C\right)^{^\frac{1}{_2}}\left(B+D\right)\left(A+C\right)^{^\frac{1}{_2}}&0\\
0&0
\end{array}\right]\right),
\end{eqnarray*}
for all $i=1,2,\cdots,2n.$ Now, we have 
\begin{eqnarray*}
&&\sum_{i=1}^{k}s_i^{2r}\left(\left[\scalebox{0.95}{$\begin{array}{lll}
B^{^\frac{1}{_2}}&0\\
D^{^\frac{1}{_2}}&0
\end{array}$}\right]\left[\begin{array}{lll}
A^{^\frac{1}{_2}}&C^{^\frac{1}{_2}}\\
0&0
\end{array}\right]\right)\\
&=&\sum_{i=1}^{k}s_i^{r}\left(\left[\begin{array}{lll}
\left(A+C\right)^{^\frac{1}{_2}}\left(B+D\right)\left(A+C\right)^{^\frac{1}{_2}}&0\\
0&0
\end{array}\right]\right)\\
&=&\sum_{i=1}^{k}s_i\left(\left[\begin{array}{lll}
\left(\left(A+C\right)^{^\frac{1}{_2}}\left(B+D\right)\left(A+C\right)^{^\frac{1}{_2}}\right)^{\frac{rp}{_p}}&0\\
0&0
\end{array}\right]\right)\\
&\leq&\sum_{i=1}^{k}s_i\left(\left[\begin{array}{lll}
\left(\left(A+C\right)^{^\frac{rp}{_2}}\left(B+D\right)^{{rp}}\left(A+C\right)^{^\frac{rp}{_2}}\right)^{\frac{1}{_p}}&0\\
0&0
\end{array}\right]\right).\\
&&\hspace{1cm}\text{(by  Lemma \ref{A417}, Lemma \ref{A409} and because $rp\geq1$)}
\end{eqnarray*}This shows that \begin{eqnarray*}
\sum_{i=1}^{k}s_i^{2r}\left(\left[\scalebox{0.95}{$\begin{array}{lll}
B^{^\frac{1}{_2}}&0\\
D^{^\frac{1}{_2}}&0
\end{array}$}\right]\left[\begin{array}{lll}
A^{^\frac{1}{_2}}&C^{^\frac{1}{_2}}\\
0&0
\end{array}\right]\right)\leq\sum_{i=1}^{k}s_i\left(\left[\scalebox{0.88}{$\begin{array}{lll}
\left(\left(A+C\right)^{^\frac{rp}{_2}}\left(B+D\right)^{{rp}}\left(A+C\right)^{^\frac{rp}{_2}}\right)^{\frac{1}{_p}}&0\\
~~0&0
\end{array}$}\right]\right),
\end{eqnarray*} 
for all $p,r\geq1$ and for all $k=1,2,\cdots,2n.$  This completes the proof.\BOX

To prove the following lemma, we need  Lemma \ref{A295}. 
\mylem{C50}{Let $A,B,C,D\in\PP_n$ and  let $r\geq0.$ Then \begin{eqnarray*}
\sum_{i=1}^{k} s_i^{2r}\left(\left[\begin{array}{lll}
A^{^\frac{1}{_2}}UB^{^\frac{1}{_2}}&0\\
0&C^{^\frac{1}{_2}}WD^{^\frac{1}{_2}}
\end{array}\right]\right)\leq\sum_{i=1}^{k}s_i^{r}\left(
\left[\begin{array}{lll}
AUB&0\\
0&CWD
\end{array}\right]\right),
\end{eqnarray*}for all $k=1,2,\cdots,2n$ and for   all unitaries $U,W\in~\MM_n(\CC)$.}
{\it Proof.~}Let $A,B,C,D\in\PP_n$ and  let $U,W$ be unitaries.  Assume that  $k=1,\cdots,2n$ and $r\geq0$. We are to prove that   \begin{eqnarray*}
\sum_{i=1}^{k} s_i^{2r}\left(\left[\begin{array}{lll}
A^{^\frac{1}{_2}}UB^{^\frac{1}{_2}}&0\\
0&C^{^\frac{1}{_2}}WD^{^\frac{1}{_2}}
\end{array}\right]\right)\leq\sum_{i=1}^{k}s_i^{r}\left(
\left[\begin{array}{lll}
AUB&0\\
0&CWD
\end{array}\right]\right).
\end{eqnarray*}
To see this, note that
\begin{eqnarray*}
&&\sum_{i=1}^{k} s_i^{2r}\left(\left[\begin{array}{lll}
A^{^\frac{1}{_2}}UB^{^\frac{1}{_2}}&0\\
0&C^{^\frac{1}{_2}}WD^{^\frac{1}{_2}}
\end{array}\right]\right)\\
&=&\sum_{i=1}^{k}\lambda_i^r\left(\left[\begin{array}{lll}
A^{^\frac{1}{_2}}UB^{^\frac{1}{_2}}&0\\
0&C^{^\frac{1}{_2}}WD^{^\frac{1}{_2}}
\end{array}\right]^*\left[\begin{array}{lll}
A^{^\frac{1}{_2}}UB^{^\frac{1}{_2}}&0\\
0&C^{^\frac{1}{_2}}WD^{^\frac{1}{_2}}
\end{array}\right]\right)\\
&=&\sum_{i=1}^{k}\lambda_i^r\left(\left[\begin{array}{lll}
B^{^\frac{1}{_2}}U^*A^{^\frac{1}{_2}}&0\\
0&D^{^\frac{1}{_2}}W^*C^{^\frac{1}{_2}}
\end{array}\right]\left[\begin{array}{lll}
A^{^\frac{1}{_2}}UB^{^\frac{1}{_2}}&0\\
0&C^{^\frac{1}{_2}}WD^{^\frac{1}{_2}}
\end{array}\right]\right)\\
&=&\sum_{i=1}^{k}\lambda_i^r\left(\left[\begin{array}{lll}
B^{^\frac{1}{_2}}&0\\
0&D^{^\frac{1}{_2}}
\end{array}\right]\left[\begin{array}{lll}
U^*&0\\
0&W^*
\end{array}\right]
\left[\begin{array}{lll}
AUB^{^\frac{1}{_2}}&0\\
0&CWD^{^\frac{1}{_2}}
\end{array}\right]\right) \\
&=&\sum_{i=1}^{k}\lambda_i^r\left(\left[\begin{array}{lll}
U^*&0\\
0&W^*
\end{array}\right]
\left[\begin{array}{lll}
AUB&0\\
0&CWD
\end{array}\right]\right) \\
&=&\sum_{i=1}^{k}\left|\lambda_i\left(\left[\begin{array}{lll}
U^*&0\\
0&W^*
\end{array}\right]
\left[\begin{array}{lll}
AUB&0\\
0&CWD
\end{array}\right]\right)\right|^r\\
&\leq&\sum_{i=1}^{k}\lambda_i^r\left(\left|\left[\begin{array}{lll}
U^*&0\\
0&W^*
\end{array}\right]
\left[\begin{array}{lll}
AUB&0\\
0&CWD
\end{array}\right]\right|\right)\hspace{0.5cm}\text{(by  Lemma  \ref{A294})}\\
&=&\sum_{i=1}^{k}s_i^r\left(\left[\begin{array}{lll}
U^*&0\\
0&W^*
\end{array}\right]
\left[\begin{array}{lll}
AUB&0\\
0&CWD
\end{array}\right]\right)\\
&=&\sum_{i=1}^{k}s_i^r\left(
\left[\begin{array}{lll}
AUB&0\\
0&CWD
\end{array}\right]\right).\hspace{0.5cm}\text{(by  Lemma \ref{A295})}\hspace{3cm}\hspace{3cm}\hspace{3cm}
\end{eqnarray*}
This shows that    \begin{eqnarray*}
\sum_{i=1}^{k} s_i^{2r}\left(\left[\begin{array}{lll}
A^{^\frac{1}{_2}}UB^{^\frac{1}{_2}}&0\\
0&C^{^\frac{1}{_2}}WD^{^\frac{1}{_2}}
\end{array}\right]\right)\leq\sum_{i=1}^{k}s_i^r\left(
\left[\begin{array}{lll}
AUB&0\\
0&CWD
\end{array}\right]\right),
\end{eqnarray*}for all $k=1,2,\cdots,2n$.  This completes the proof.\BOX

 Let us now prove our  main result  using
Lemma \ref{C51} and\text{ Lemma} \ref{C50}.  Recall that the Cauchy-Schwarz inequality states that 
\begin{eqnarray}\label{C25}
\sum_{i=1}^n|x_iy_i|\leq\left(\sum_{i=1}^nx_i^2\right)^\frac{1}{2}\left(\sum_{i=1}^ny_i^2\right)^\frac{1}{2} \text{ for all $x_i,y_i\in\RR,~i=1,\cdots,n.$}
\end{eqnarray}
\mythe{ASAS1}{(The Main Theorem) Let  $A,B,C,D\in\PP_n.$  Then     \begin{eqnarray*}
\hspace{0.5cm} s\left(\left(A^{^2}\sharp B^{^2}\right)^{r}+\left(C^{^2}\sharp D^{^2}\right)^{r}\right)\prec_w s\left(\left(\left(A+C\right)^{^\frac{rp}{_2}}\left(B+D\right)^{{rp}}\left(A+C\right)^{^\frac{rp}{_2}}\right)^{\frac{1}{_p}}\right),
\end{eqnarray*} for all $p\geq1$ and for all $r\geq1$.}
{\it Proof.~} Let  $A,B,C,D\in\PP_n$ and let $p\geq1$ and $r\geq1.$   Then, using Lemma \ref{C2}, we have $A^{^2}\sharp B^{^2}=AUB$ and $C^{^2}\sharp D^{^2}=CWD$ for some unitaries  $U,W.$     We define  $$X=\left(A^{^2}\sharp B^{^2}\right)^{r}+\left(C^{^2}\sharp D^{^2}\right)^{r}$$ and $$Y=\left(\left(A+C\right)^{^\frac{rp}{_2}}\left(B+D\right)^{{rp}}\left(A+C\right)^{^\frac{rp}{_2}}\right)^{\frac{1}{_p}}.$$ We are to prove that $s\left(X\right)\prec_w s\left(Y\right)$. For all~$k=1,2,\cdots,2n,$ we have
\begin{eqnarray*}
&&\sum_{i=1}^{k} s_i\left(\left[\begin{array}{lll}
X&0\\
0&0
\end{array}\right]\right)\\
&=&\sum_{i=1}^{k} s_i\left(\left[\begin{array}{lll}
\left(A^{^2}\sharp B^{^2}\right)^{r}+\left(C^{^2}\sharp D^{^2}\right)^{r}&0\\
0&0
\end{array}\right]\right)\\
&\leq&\sum_{i=1}^{k} s_i\left(\left[\begin{array}{lll}
\left(A^{^2}\sharp B^{^2}+C^{^2}\sharp D^{^2}\right)^{r}&0\\
0&0
\end{array}\right]\right)\\
&&\text{(because $f(t)=t^r$ is convex  and by Lemma \ref{A389} and  Lemma \ref{C32})}\\
&=&\sum_{i=1}^{k} s_i^{r}\left(\left[\begin{array}{lll}
A^{^2}\sharp B^{^2}+C^{^2}\sharp D^{^2}&0\\
0&0
\end{array}\right]\right)\\
&=&\sum_{i=1}^{k}s_i^{r}\left(
\left[\begin{array}{lll}
AUB+CWD&0\\
0&0
\end{array}\right]\right)\\
&=&\sum_{i=1}^{k}s_i^{r}\left(
\left[\begin{array}{lll}
A^{^\frac{1}{_2}}&C^{^\frac{1}{_2}}\\
0&0
\end{array}\right]\left[\begin{array}{lll}
A^{^\frac{1}{_2}}UB^{^\frac{1}{_2}}&0\\
0&C^{^\frac{1}{_2}}WD^{^\frac{1}{_2}}
\end{array}\right]\left[\begin{array}{lll}
B^{^\frac{1}{_2}}&0\\
D^{^\frac{1}{_2}}&0
\end{array}\right]\right)\\
&\leq&\sum_{i=1}^{k}s_i^{r}\left(
\left[\begin{array}{lll}
A^{^\frac{1}{_2}}UB^{^\frac{1}{_2}}&0\\
0&C^{^\frac{1}{_2}}WD^{^\frac{1}{_2}}
\end{array}\right]\left[\begin{array}{lll}
B^{^\frac{1}{_2}}&0\\
D^{^\frac{1}{_2}}&0
\end{array}\right]\left[\begin{array}{lll}
A^{^\frac{1}{_2}}&C^{^\frac{1}{_2}}\\
0&0
\end{array}\right]\right)~\\
&&\hspace{3cm}\text{(by Lemma \ref{A306},  Lemma \ref{C32}  and Lemma  \ref{A341})}\\
&\leq&\sum_{i=1}^{k}s_i^{r}\left(
\left[\begin{array}{lll}
A^{^\frac{1}{_2}}UB^{^\frac{1}{_2}}&0\\
0&C^{^\frac{1}{_2}}WD^{^\frac{1}{_2}}
\end{array}\right]\right)s_i^{r}\left(\left[\begin{array}{lll}
B^{^\frac{1}{_2}}&0\\
D^{^\frac{1}{_2}}&0
\end{array}\right]\left[\begin{array}{lll}
A^{^\frac{1}{_2}}&C^{^\frac{1}{_2}}\\
0&0
\end{array}\right]\right)\\
&&\hspace{8.4cm}\text{(by Lemma \ref{A298})}\hspace{9cm}
\end{eqnarray*}
\begin{eqnarray*}
&=&\left(\sum_{i=1}^{k}s_i^{2r}\left(
\left[\scalebox{0.95}{$\begin{array}{lll}
A^{^\frac{1}{_2}}UB^{^\frac{1}{_2}}&0\\
0&C^{^\frac{1}{_2}}WD^{^\frac{1}{_2}}
\end{array}$}\right]\right)\right)^\frac{1}{_2}\left(\sum_{i=1}^{k}s_i^{2r}\left(\left[\scalebox{0.95}{$\begin{array}{lll}
B^{^\frac{1}{_2}}&0\\
D^{^\frac{1}{_2}}&0
\end{array}$}\right]\left[\begin{array}{lll}
A^{^\frac{1}{_2}}&C^{^\frac{1}{_2}}\\
0&0
\end{array}\right]\right)\right)^\frac{1}{_2}\\
&&\hspace{10cm}\text{(by (\ref{C25}))}\\
&\leq&\left(\sum_{i=1}^{k}s_i^{2r}\left(
\left[\scalebox{0.75}{$\begin{array}{lll}
A^{^\frac{1}{_2}}UB^{^\frac{1}{_2}}&0\\
0&C^{^\frac{1}{_2}}WD^{^\frac{1}{_2}}
\end{array}$}\right]\right)\right)^\frac{1}{_2}\left(\sum_{i=1}^{k}s_i\left(\scalebox{0.7}{$\left[\begin{array}{lll}
\scalebox{1.1}{$\left(\left(A+C\right)^{\scalebox{1}{$^\frac{rp}{_2}$}}\left(B+D\right)^{\scalebox{0.9}{$rp$}}\left(A+C\right)^{\scalebox{1}{$^\frac{rp}{_2}$}}\right)^{\frac{1}{_p}}$}&0\\
0&0
\end{array}\right]$}\right)\right)^\frac{1}{_2}\\
&&\hspace{9cm}\text{(by Lemma \ref{C51})}\\
&=&\left(\sum_{i=1}^{k}s_i^{2r}\left(
\left[\begin{array}{lll}
A^{^\frac{1}{_2}}UB^{^\frac{1}{_2}}&0\\
0&C^{^\frac{1}{_2}}WD^{^\frac{1}{_2}}
\end{array}\right]\right)\right)^\frac{1}{_2}\left(\sum_{i=1}^{k}s_i\left(\left[\begin{array}{lll}
Y&0\\
~~0&0
\end{array}\right]\right)\right)^\frac{1}{_2}\\
&\leq&\left(\sum_{i=1}^{k}s_i^{r}\left(
\left[\begin{array}{lll}
AUB&0\\
0&CWD
\end{array}\right]\right)\right)^\frac{1}{_2}\left(\sum_{i=1}^{k}s_i\left(\left[\begin{array}{lll}
Y&0\\
~~0&0
\end{array}\right]\right)\right)^\frac{1}{_2}\hspace{0.5cm}\text{(by Lemma \ref{C50})}\\
&=&\left(\sum_{i=1}^{k}s_i^{r}\left(
\left[\begin{array}{lll}
A^{^2}\sharp B^{^2}&0\\
0&C^{^2}\sharp D^{^2}
\end{array}\right]\right)\right)^\frac{1}{_2}\left(\sum_{i=1}^{k}s_i\left(\left[\begin{array}{lll}
Y&0\\
~~0&0
\end{array}\right]\right)\right)^\frac{1}{_2}\\
&=&\left(\sum_{i=1}^{k}s_i\left(
\left[\begin{array}{lll}
\left(A^{^2}\sharp B^{^2}\right)^{r}&0\\
0&\left(C^{^2}\sharp D^{^2}\right)^{r}
\end{array}\right]\right)\right)^\frac{1}{_2}\left(\sum_{i=1}^{k}s_i\left(\left[\begin{array}{lll}
Y&0\\
~~0&0
\end{array}\right]\right)\right)^\frac{1}{_2}\\
&\leq&\left(\sum_{i=1}^{k}s_i\left(
\left[\begin{array}{lll}
\left(A^{^2}\sharp B^{^2}\right)^{r}+\left(C^{^2}\sharp D^{^2}\right)^{r}&0\\
0&0
\end{array}\right]\right)\right)^\frac{1}{_2}\left(\sum_{i=1}^{k}s_i\left(\left[\begin{array}{lll}
Y&0\\
~~0&0
\end{array}\right]\right)\right)^\frac{1}{_2}\\
&&\hspace{7cm}\text{(by Lemma \ref{A299})}\\
&=&\left(\sum_{i=1}^{k}s_i\left(
\left[\begin{array}{lll}
X&0\\
0&0
\end{array}\right]\right)\right)^\frac{1}{_2}\left(\sum_{i=1}^{k}s_i\left(\left[\begin{array}{lll}
Y&0\\
~~0&0
\end{array}\right]\right)\right)^\frac{1}{_2}.\hspace{9cm}
\end{eqnarray*}
This shows that \begin{eqnarray*}
\displaystyle\sum_{i=1}^{k} s_i\left(\left[\begin{array}{lll}
X&0\\
0&0
\end{array}\right]\right)\leq\left(\displaystyle\sum_{i=1}^{k} s_i\left(\left[\begin{array}{lll}
X&0\\
0&0
\end{array}\right]\right)\right)^\frac{1}{_2}\left(\displaystyle\sum_{i=1}^{k}s_i\left(\left[\begin{array}{lll}
Y&0\\
~~0&0
\end{array}\right]\right)\right)^\frac{1}{_2},
\end{eqnarray*}
for all $k=1,\cdots,2n.$
Thus we have
\begin{eqnarray*}
\sum_{i=1}^{k}s_i\left(
X\oplus0\right)\leq\left(\sum_{i=1}^{k}s_i\left(
X\oplus0\right)\right)^\frac{1}{_2}\left(\sum_{i=1}^{k}s_i\left(Y\oplus0\right)\right)^\frac{1}{_2},
\end{eqnarray*}for all $k=1,2,\cdots,2n.$
Hence\begin{eqnarray*}
\left(\sum_{i=1}^{k}s_i\left(
X\oplus0\right)\right)^\frac{1}{_2}\leq\left(\sum_{i=1}^{k}s_i\left(Y\oplus0\right)\right)^\frac{1}{_2},
\end{eqnarray*}for all $k=1,2,\cdots,2n.$
From this we get $$\sum_{i=1}^{k}s_i\left(
X\oplus0\right)\leq\sum_{i=1}^{k}s_i\left(Y\oplus0\right),$$for all $k=1,2,\cdots,2n.$
So we have $$s\left(X\oplus0\right)\prec_ws\left(Y\oplus0\right).$$
Thus $$s\left(X\right)\prec_ws\left(Y\right).$$
This is equivalent to$$s\left(\left(A^{^2}\sharp B^{^2}\right)^{r}+\left(C^{^2}\sharp D^{^2}\right)^{r}\right)\prec_ws\left(\left(\left(A+C\right)^{^\frac{rp}{_2}}\left(B+D\right)^{{rp}}\left(A+C\right)^{^\frac{rp}{_2}}\right)^{\frac{1}{_p}}\right),$$
for all $p,r\geq1$. This completes the proof.\BOX

\mycor{}{Let  $A,B,C,D\in\PP_n$  and let $p\geq1$ and $r\geq1.$ Then  $$\left|\left|\left|\left(A^{^2}\sharp B^{^2}\right)^{r}+\left(C^{^2}\sharp D^{^2}\right)^{r}\right|\right|\right|\leq \left|\left|\left|\left(\left(A+C\right)^{^\frac{rp}{_2}}\left(B+D\right)^{{rp}}\left(A+C\right)^{^\frac{rp}{_2}}\right)^{\frac{1}{_p}}\right|\right|\right|,$$for all  unitarly invariant norms.}
{\it Proof.~}By using Theorem \ref{ASAS1} and  Lemma \ref{C32}, we are done.  \BOX

We remark here that the previous corollary settles Conjecture \ref{CON5} for the case of $m=2$ and settles Conjecture \ref{CON2} for the case of $t=\frac{1}{2},$   $r\geq1,$  $ p\geq1$ and $m=2.$
\section{Trace Norm for Sum of Geometric Means}\label{S33}
We are now in a position to prove Conjecture \ref{CON2} in the case of the trace norm.  
\mythe{A496}{Let $A_i,B_i\in\PP_n$ for all $i=1,\cdots,m.$ Then $$\left|\left|\sum_{i=1}^m(A_i^2\sharp B_i^2)^r\right|\right|_1\leq\left|\left|\left(\left(\sum_{i=1}^mA_i\right)^{\frac{pr}{_2}}\left(\sum_{i=1}^mB_i\right)^{pr}\left(\sum_{i=1}^mA_i\right)^{\frac{rp}{_2}}\right)^{\frac{1}{p}}\right|\right|_1,$$for all $p>0$ and for all $r\geq1.$}
{\it Proof.~}Let $A_i,B_i\in\PP_n$ for all $i=1,\cdots,m.$ Let $r\geq1$ and $p>0$. Using\text{ Lemma}~\ref{A403}, we have \begin{eqnarray*}
\lambda\left(A_i^2\sharp B_i^2\right)\prec_{\log}\lambda\left((A_i\sharp B_i)^2\right)~~\text{for all }i=1,\cdots,m.
\end{eqnarray*}
Then, using Lemma \ref{A409}, we have \begin{eqnarray*}
\lambda\left(A_i^2\sharp B_i^2\right)\prec_w\lambda\left((A_i\sharp B_i)^2\right)~~\text{for all }i=1,\cdots,m.
\end{eqnarray*}
 Using Lemma \ref{A341},  we have
\begin{eqnarray}\label{A476}
\lambda^r\left(A_i^2\sharp B_i^2\right)\prec_w\lambda^r\left((A_i\sharp B_i)^2\right)~~\text{for all }i=1,\cdots,m.
\end{eqnarray}
Thus we have \begin{eqnarray}\label{A407}
\lambda\left(\left(A_i^2\sharp B_i^2\right)^r\right)\prec_w\lambda\left((A_i\sharp B_i)^{2r}\right)~~\text{for all }i=1,\cdots,m.
\end{eqnarray}
Hence \begin{eqnarray}\label{A481}
s\left(\left(A_i^2\sharp B_i^2\right)^r\right)\prec_ws\left((A_i\sharp B_i)^{2r}\right)~~\text{for all }i=1,\cdots,m.
\end{eqnarray}
From Lemma  \ref{C32}, we get \begin{eqnarray*}\label{A482}
\left|\left|\left|\left(A_i^2\sharp B_i^2\right)^r\right|\right|\right|\leq\left|\left|\left|(A_i\sharp B_i)^{2r}\right|\right|\right|,
\end{eqnarray*}
$\text{for all }i=1,\cdots,m\text{ and  for any unitarly invariant norms.}$ In particular, $$\left|\left|\left(A_i^2\sharp B_i^2\right)^r\right|\right|_1\leq\left|\left|(A_i\sharp B_i)^{2r}\right|\right|_1~~\text{for all }i=1,\cdots,m.$$
 Hence \begin{eqnarray}\label{A920}
 \sum_{i=1}^m\left|\left|\left(A_i^2\sharp B_i^2\right)^r\right|\right|_1\leq\sum_{i=1}^m\left|\left|(A_i\sharp B_i)^{2r}\right|\right|_1.
 \end{eqnarray}
 Using (\ref{A920}), we have 
$$\left|\left|\sum_{i=1}^m\left(A_i^2\sharp B_i^2\right)^r\right|\right|_1=\sum_{i=1}^m\left|\left|\left(A_i^2\sharp B_i^2\right)^r\right|\right|_1\leq\sum_{i=1}^m\left|\left|(A_i\sharp B_i)^{2r}\right|\right|_1=\left|\left|\sum_{i=1}^m(A_i\sharp B_i)^{2r}\right|\right|_1.$$This shows that
\begin{eqnarray}\label{A491}
\left|\left|\sum_{i=1}^m\left(A_i^2\sharp B_i^2\right)^r\right|\right|_1\leq\left|\left|\sum_{i=1}^m(A_i\sharp B_i)^{2r}\right|\right|_1.
\end{eqnarray}
Now, we have 
\begin{eqnarray*}
\left|\left|\sum_{i=1}^m(A_i^2\sharp B_i^2)^r\right|\right|_1&\leq&\left|\left|\sum_{i=1}^m(A_i\sharp B_i)^{2r}\right|\right|_1\hspace{0.5cm}\text{(by (\ref{A491}))}\\
&\leq&\left|\left|\left(\left(\sum_{i=1}^mA_i\right)^{\frac{pr}{_2}}\left(\sum_{i=1}^mB_i\right)^{pr}\left(\sum_{i=1}^mA_i\right)^{\frac{rp}{_2}}\right)^{\frac{1}{p}}\right|\right|_1.\\
&&\text{(by taking  $t=\frac{1}{2}$ and replacing $r$ by $2r$ in inequality (\ref{A2}))}
\end{eqnarray*}
This shows that $$\left|\left|\sum_{i=1}^m(A_i^2\sharp B_i^2)^r\right|\right|_1\leq\left|\left|\left(\left(\sum_{i=1}^mA_i\right)^{\frac{pr}{_2}}\left(\sum_{i=1}^mB_i\right)^{pr}\left(\sum_{i=1}^mA_i\right)^{\frac{rp}{_2}}\right)^{\frac{1}{p}}\right|\right|_1,$$
for all $p>0$ and for all $r\geq1.$  This completes the proof.\BOX

\section{A New Proof of Some Special Cases of  Inequality (\ref{A2})}\label{S44}
In this section, we  present  a new proof  of  a special case  of  inequality (\ref{A2}) when $p=m=2,$ $t=\frac{1}{2},$ and  $r\geq1$.  First, we give a new proof of   inequality (\ref{A2})  when $p=m=2,$  $t=\frac{1}{2}$  and  $r=1$ as in the following theorem.
\mythe{A340}{Let  $A,B,C,D\in\PP_n.$ Then \begin{eqnarray}
\left|\left|\left|A\sharp B+C\sharp D\right|\right|\right|\leq\left|\left|\left|\left(\left(A+C\right)^{^\frac{1}{_2}}\left(B+D\right)\left(A+C\right)^{^\frac{1}{_2}}\right)^\frac{1}{2}\right|\right|\right|,
\end{eqnarray}for all  unitarly invariant norms.}
{\it Proof.~}Let $A,B,C,D\in\PP_n.$  Using Lemma \ref{C2}, we have  $A\sharp B=A^{^\frac{1}{_2}}UB^{^\frac{1}{_2}}$ and $C\sharp D=C^{^\frac{1}{_2}}WD^{^\frac{1}{_2}}$ for some unitaries  $U,W\in~\MM_n(\CC).$ We are to prove that $\left|\left|\left|A\sharp B+C\sharp D\right|\right|\right|\leq\left|\left|\left|\left(\left(A+C\right)^{^\frac{1}{_2}}\left(B+D\right)\left(A+C\right)^{^\frac{1}{_2}}\right)^\frac{1}{2}\right|\right|\right|$ for all  unitarly invariant norms. Then for all $k=1,\cdots,2n,$ we have
 \begin{eqnarray*}
 &&\sum_{i=1}^{k} s_i\left(\left[\begin{array}{lll}
A\sharp B+C\sharp D&0\\
0&0
\end{array}\right]\right)\\
&=&\sum_{i=1}^{k} s_i\left(\left[\begin{array}{lll}
A^{^\frac{1}{_2}}UB^{^\frac{1}{_2}}+C^{^\frac{1}{_2}}WD^{^\frac{1}{_2}}&0\\
0&0
\end{array}\right]\right)\\
&=&\sum_{i=1}^{k} s_i\left(\left[\begin{array}{lll}
A^{^\frac{1}{_2}}U&C^{^\frac{1}{_2}}W\\
0&0
\end{array}\right]\left[\begin{array}{lll}
B^{^\frac{1}{_2}}&0\\
D^{^\frac{1}{_2}}&0
\end{array}\right]\right)\\
&=&\sum_{i=1}^{k} s_i\left(\left[\begin{array}{lll}
A^{^\frac{1}{_2}}&C^{^\frac{1}{_2}}\\
0&0
\end{array}\right]\left[\begin{array}{lll}
U&0\\
0&W
\end{array}\right]\left[\begin{array}{lll}
B^{^\frac{1}{_2}}&0\\
D^{^\frac{1}{_2}}&0
\end{array}\right]\right)\\
&\leq&\sum_{i=1}^{k} s_i\left(\left[\begin{array}{lll}
U&0\\
0&W
\end{array}\right]\left[\begin{array}{lll}
B^{^\frac{1}{_2}}&0\\
D^{^\frac{1}{_2}}&0
\end{array}\right]\left[\begin{array}{lll}
A^{^\frac{1}{_2}}&C^{^\frac{1}{_2}}\\
0&0
\end{array}\right]\right)\\
&&\hspace{0.5cm}\hspace{0.5cm}\hspace{0.5cm}\hspace{0.5cm}\hspace{0.5cm}\hspace{0.5cm}\text{(by Lemma \ref{C32} and  Lemma \ref{A306})}\\
&=&\sum_{i=1}^{k} s_i\left(\left[\begin{array}{lll}
B^{^\frac{1}{_2}}&0\\
D^{^\frac{1}{_2}}&0
\end{array}\right]\left[\begin{array}{lll}
A^{^\frac{1}{_2}}&C^{^\frac{1}{_2}}\\
0&0
\end{array}\right]\right).\\
&&\hspace{1.5cm}\text{\Big(by Lemma \ref{A295} and because  $\left[\begin{array}{lll}
U&0\\
0&W
\end{array}\right]$ is unitary\Big)}\\
&=&\sum_{i=1}^{k} \left(s_i^2\left(\left[\begin{array}{lll}
B^{^\frac{1}{_2}}&0\\
D^{^\frac{1}{_2}}&0
\end{array}\right]\left[\begin{array}{lll}
A^{^\frac{1}{_2}}&C^{^\frac{1}{_2}}\\
0&0
\end{array}\right]\right)\right)^\frac{1}{2}\\
 &=&\sum_{i=1}^k\left(s_i\left(\left[\begin{array}{lll}
\left(A+C\right)^{^\frac{1}{_2}}\left(B+D\right)\left(A+C\right)^{^\frac{1}{_2}}&0\\
~~0&0
\end{array}\right]\right)\right)^\frac{1}{_2}\hspace{0.5cm}\text{(by Lemma \ref{C53})}\\
&=&\sum_{i=1}^ks_i\left(\left[\begin{array}{lll}
\left(\left(A+C\right)^{^\frac{1}{_2}}\left(B+D\right)\left(A+C\right)^{^\frac{1}{_2}}\right)^\frac{1}{_2}&0\\
~~0&0
\end{array}\right]\right).
 \end{eqnarray*}
 From this, we get   
  \begin{eqnarray*}
\sum_{i=1}^{k} s_i\left(\left[\begin{array}{lll}
A\sharp B+C\sharp D&0\\
0&0
\end{array}\right]\right)\leq\sum_{i=1}^ks_i\left(\left[\begin{array}{lll}
\left(\left(A+C\right)^{^\frac{1}{_2}}\left(B+D\right)\left(A+C\right)^{^\frac{1}{_2}}\right)^\frac{1}{2}&0\\
~~0&0
\end{array}\right]\right),
 \end{eqnarray*}
for all $k=1,\cdots,2n.$ Thus we have 
 $$s\left((A\sharp B+C\sharp D)\oplus0\right)\prec_ws\left(\left(\left(A+C\right)^{^\frac{1}{_2}}\left(B+D\right)\left(A+C\right)^{^\frac{1}{_2}}\right)^\frac{1}{2}\oplus0\right).$$
Hence\begin{eqnarray*}
s\left(A\sharp B+C\sharp D\right)\prec_ws\left(\left(\left(A+C\right)^{^\frac{1}{_2}}\left(B+D\right)\left(A+C\right)^{^\frac{1}{_2}}\right)^\frac{1}{2}\right)
\end{eqnarray*} 
Using Lemma \ref{C32}, we have 
 $$\left|\left|\left|A\sharp B+C\sharp D\right|\right|\right|\leq\left|\left|\left|\left(\left(A+C\right)^{^\frac{1}{_2}}\left(B+D\right)\left(A+C\right)^{^\frac{1}{_2}}\right)^\frac{1}{2}\right|\right|\right|,$$
for all  unitarly invariant norms. This completes the proof.\BOX

Second, we give a new proof of   inequality (\ref{A2})  when $p=m=2,$ $t=\frac{1}{2}$ and  $r\geq1$   using Theorem \ref{A340}. 
  \mythe{A495}{Let $A,B,C,D\in\PP_n$ and $r\geq1.$ Then $$\left|\left|\left|(A\sharp B)^r+(C\sharp D)^r\right|\right|\right|\leq\left|\left|\left|\left(\left(A+C\right)^\frac{r}{_2}\left(B+D\right)^r\left(A+C\right)^{^\frac{r}{_2}}\right)^{^\frac{1}{_2}}\right|\right|\right|,$$for all  unitarly invariant norms.}
{\it Proof.~}Let $A,B,C,D\in\PP_n$ and $r\geq1.$  We are to prove that $$\left|\left|\left|(A\sharp B)^r+(C\sharp D)^r\right|\right|\right|\leq\left|\left|\left|\left(\left(A+C\right)^\frac{r}{_2}\left(B+D\right)^r\left(A+C\right)^{^\frac{r}{_2}}\right)^{^\frac{1}{_2}}\right|\right|\right|,$$for all  unitarly invariant norms.  
Using  Lemma \ref{A389} and since  the function~$f:[0,\infty)\longrightarrow[0,\infty)$ given by $f(x)=x^r$   is convex, it follows that
\begin{eqnarray}\label{A922}
\left|\left|\left|(A\sharp B)^r+(C\sharp D)^r\right|\right|\right|\leq\left|\left|\left|(A\sharp B+C\sharp D)^r\right|\right|\right|,\end{eqnarray}for all  unitarly invariant norms. Using Theorem \ref{A340}, we have
  $$\left|\left|\left|A\sharp B+C\sharp D\right|\right|\right|\leq\left|\left|\left|\left(\left(A+C\right)^{^\frac{1}{_2}}\left(B+D\right)\left(A+C\right)^{^\frac{1}{_2}}\right)^\frac{1}{2}\right|\right|\right|,$$for all  unitarly invariant norms. 
From this, we note that
\begin{eqnarray*}
&&\left|\left|\left|A\sharp B+C\sharp D\right|\right|\right|\leq\left|\left|\left|\left(\left(A+C\right)^{^\frac{1}{_2}}\left(B+D\right)\left(A+C\right)^{^\frac{1}{_2}}\right)^\frac{1}{2}\right|\right|\right|\\
&\Longrightarrow&\left|\left|\left|\left( A\sharp B+C\sharp D\right)^r\right|\right|\right|\leq\left|\left|\left|\left(\left(\left(A+C\right)^{^\frac{1}{_2}}\left(B+D\right)\left(A+C\right)^{^\frac{1}{_2}}\right)^\frac{1}{2}\right)^r\right|\right|\right|\\
&&\hspace{8cm}\text{(by Lemma \ref{A429})}\\
&\Longrightarrow&\left|\left|\left|\left( A\sharp B+C\sharp D\right)^r\right|\right|\right|\leq\left|\left|\left|\left(\left(A+C\right)^{^\frac{1}{_2}}\left(B+D\right)\left(A+C\right)^{^\frac{1}{_2}}\right)^\frac{r}{_2}\right|\right|\right|.
\end{eqnarray*}
So \begin{eqnarray}\label{A924}
\left|\left|\left|\left( A\sharp B+C\sharp D\right)^r\right|\right|\right|\leq\left|\left|\left|\left(\left(A+C\right)^{^\frac{1}{_2}}\left(B+D\right)\left(A+C\right)^{^\frac{1}{_2}}\right)^\frac{r}{_2}\right|\right|\right|,
\end{eqnarray}
for all  unitarly invariant norms. From Lemma \ref{A417}, Lemma \ref{A409} and Lemma \ref{C32}, we get \begin{eqnarray}\label{Ass924}
\left|\left|\left|\left(\left(A+C\right)^{^\frac{1}{_2}}\left(B+D\right)\left(A+C\right)^{^\frac{1}{_2}}\right)^\frac{r}{_2}\right|\right|\right|\leq\left|\left|\left|\left(\left(A+C\right)^{^\frac{r}{_2}}\left(B+D\right)^r\left(A+C\right)^{^\frac{r}{_2}}\right)^\frac{1}{_2}\right|\right|\right|,\hspace{0.5cm}
\end{eqnarray}
for all  unitarly invariant norms. 
  Now, using (\ref{A924}), (\ref{A922}) and (\ref{Ass924}), we conclude that 
 $$\left|\left|\left|(A\sharp B)^r+(C\sharp D)^r\right|\right|\right|\leq\left|\left|\left|\left(\left(A+C\right)^{^\frac{r}{_2}}\left(B+D\right)^r\left(A+C\right)^{^\frac{r}{_2}}\right)^\frac{1}{_2}\right|\right|\right|,$$
for all  unitarly invariant norms and for $r\geq1$. This completes the proof.\BOX

\end{document}